\documentclass[12pt]{amsart}
\numberwithin{equation}{section}
\newtheorem{theorem}{Theorem}
\numberwithin{theorem}{section}

\numberwithin{theorem}{section} \numberwithin{lemma}{section}

\numberwithin{definition}{section}
\newtheorem{corollary}{Corollary}
\numberwithin{corollary}{section}

\numberwithin{remark}{section}

\numberwithin{proposition}{section}
\def\b{\begin{equation}}
\def\e{\end{equation}}
\newcommand{\ignore}[1]{}
\normalsize
\date {November 27, 2006}
\begin{document}
\pagenumbering{arabic} \pagenumbering{arabic}\setcounter{page}{1}
\tracingpages 1
\title{Hardy, Rellich  and Uncertainty principle inequalities  on Carnot Groups}
\author{Ismail Kombe }
\dedicatory {}
\address{Ismail Kombe, Mathematics Department\\ Dawson-Loeffler Science
\&Mathematics Bldg\\
Oklahoma City University \\
2501 N. Blackwelder, Oklahoma City, OK 73106-1493}
\email{ikombe@okcu.edu}
\begin{abstract}
In this paper we prove sharp weighted Hardy-type inequalities on
Carnot groups with the homogeneous norm $N=u^{1/(2-Q)}$ associated
to Folland's fundamental solution $u$ for the sub-Laplacian
$\Delta_{\mathbb{G}}$. We also prove uncertainty principle,
Caffarelli-Kohn-Nirenberg and Rellich inequalities on Carnot
groups.
\end{abstract}
\maketitle
\section{Introduction}

The classical Hardy inequality states that for $n\ge 3$
\begin{equation}
\int_{\mathbb{R}^n}|\nabla\phi(x)|^2dx\ge
\Big(\frac{n-2}{2}\Big)^2\int_{\mathbb{R}^n}
\frac{|\phi(x)|^2}{|x|^2}dx,
\end{equation}
where $\phi\in C_0^{\infty}( \mathbb{R}^n\setminus \{0\})$ and the
constant $(\frac{n-2}{2})^2$ is sharp. There exists a large
literature dealing with the Hardy-type inequalities on the
Euclidean space $\mathbb{R}^n$ and, in particular, sharp
inequalities as well as their improved versions which have
attracted a lot of  attention because of their application to
singular problems, e.g.  \cite{4}, \cite{8}, \cite{9}, \cite{37},
\cite{42}. For instance, Baras and Goldstein in their classical
paper \cite {4}, showed that the following heat problem
\begin{equation}
\begin{cases}
\frac{\partial u}{\partial t}=\Delta u +\frac{c}{|x|^2}u \quad
&\text{in}\quad \Omega\times(0, \infty), \quad 0\in \Omega, \\
u(x,0)=u_0(x)\ge 0 \quad &\text{in}  \quad \Omega,\\
u(x,t)=0 \quad &\text{on} \quad \partial \Omega\times (0,\infty),
\end{cases}
\end{equation}
has a global solution ( in the sense of distributions) if $c \le
C^*(n)=(\frac{n-2}{2})^2$ and no solution, even locally in time,
if $ c>C^*(n)=(\frac{n-2}{2})^2$. Thus, $C^*(n)=(\frac{n-2}{2})^2$
is the cut-off point for existence of positive solutions for the
heat equation with inverse square potential $c/|x|^2$.

Recently there has been considerable interest in improving
 the inequality (1.1), in the sense that nonnegative terms are added in the right hand side of
 (1.1), and one of the important improvement has been obtained by  Brezis and
 V\'azquez \cite {8}. They proved that for a bounded domain $\Omega\subset
 \mathbb{R}^n$ there holds
 \begin{equation}
\int _{\Omega} |\nabla \phi(x)|^2dx\ge
\Big(\frac{n-2}{2}\Big)^2\int_{\Omega}
\frac{|\phi(x)|^2}{|x|^2}dx+\mu\big(\frac{\omega_n}{|\Omega|}\big)^{2/n}\int
_{\Omega} \phi^2dx,
 \end{equation}
where $\omega_n$  and $|\Omega|$ denote the $n$-dimensional
Lebesgue measure of the unit ball $B\subset \mathbb{R}^n$  and the
domain $\Omega$ respectively. Here $\mu=5.7832$ is the first
eigenvalue of the Laplace operator in the two dimensional unit
disk and it is optimal when $\Omega$ is a ball centered at the
origin. A comprehensive treatment of improved Hardy inequalities
with best constants, involving various kinds of distance functions
in the Euclidean space $ \mathbb{R}^n$  can be found in \cite {6}.

In view of  these important works mentioned above it is natural to
investigate Hardy-type inequalities and their improved versions on
general Carnot groups. It is well known that the Euclidean space $
\mathbb{R}^n$ with its usual abelian group structure is a trivial
Carnot group. We are mainly concerned with the Hardy-type
inequalities on non-trivial Carnot groups.

The simplest nontrivial example of a Carnot group is given by the
Heisenberg group $ \mathbb{H}^n$. The following Hardy-type
inequality on the Heisenberg group $ \mathbb{H}^n$ was first
proved by Garofalo and Lanconelli \cite {23} :
\begin{equation}
\int_{ \mathbb{H}^n}|\nabla_{ \mathbb{H}^n}\phi|^2dzdt \ge
\Big(\frac{Q-2}{2}\Big)^2\int_{ \mathbb{H}^n}
(\frac{|z|^2}{|z|^4+t^2})\phi^2 dzdt
\end{equation}
where $\phi\in C_0^{\infty}(\mathbb{H}^n\setminus \{0\})$,
$Q=2n+2$ and the constant $(\frac{Q-2}{2})^2$ is sharp. Here we
view $\mathbb{H}^n$ as $ \mathbb{C}^n\times\mathbb{R}$, and $dzdt$
refers to the usual Lebesgue measure.  Further results concerning
Hardy-type inequality on the Heisenberg group can be found in
\cite {36} and \cite {13}. Recently, Han and Niu \cite {26}, and
D'Ambrosio \cite {14} obtained a version of Hardy-Sobolev
inequality on the \textit{H}-type group  and Hardy-type
inequalities on Carnot groups, respectively. We indicate that a
result in  \cite {14} concerning Hardy-type inequality  on general
Carnot groups overlap with ours (Theorem 4.1), but the methods of
proof are different.

The first goal of this paper is to investigate the existence and
the explicit determination of constants $C$ and weight $q(x)$ on
Carnot group $\mathbb{G}$ such that the Hardy-type inequality

\begin{equation}
\int_{ \mathbb{G}}w(x)|\nabla_{ \mathbb{G}} \phi(x)|^2dx \ge
C\int_{ \mathbb{G}} q(x)|\phi(x)|^2dx
\end{equation}
holds for all $\phi\in C_0^{\infty}(\mathbb{G} \setminus \{0\})$.
Here we consider a special weight function $w(x)$ which is related
to the fundamental solution of sub-Laplacian $ \Delta_{
\mathbb{G}}$ on Carnot group $\mathbb{G}$, and $dx$ refers to Haar
measure on $\mathbb{G}$.

It is important to emphasize that our result lead us to obtain a
version of the uncertainty principle, Caffarelli-Kohn-Nirenberg
and Rellich inequalities on general Carnot groups.

Although we prove Hardy-type inequalities on Carnot groups with an
arbitrary step, we first establish sharp  Hardy-type inequalities
on the Heisenberg group $ \mathbb{H}^n$ and extend this result to
the \textit{H}-type groups. The main reason for doing this is that
the fundamental solution of the sub-Laplacian on the Heisenberg
group $ \mathbb{H}^n$  and \textit{H}-type groups are known
explicitly (see Section 3) but not for general Carnot groups. The
proof of our theorem on general Carnot groups differs slightly in
some steps from the Heisenberg group $ \mathbb{H}^n$ and
\textit{H}-type group cases. The method that we apply here,
inspired by the work of Allegretto \cite {2}, can be applied to
the Baouendi-Grushin type vector fields in that they do not arise
from any Carnot group.

The plan of the paper is as follows: In Section 2, we recall the
basic properties of Carnot group $\mathbb{G}$ and some well known
results that will be used in the sequel. In Section 3, we prove
the Hardy-type inequalities on the Heisenberg group and
\textit{H}-type group. In Section 4, we prove Hardy-type
inequality on general Carnot groups. As a consequence of the
Hardy-type inequality, we obtain a version of uncertainty
principle and Caffarelli-Kohn-Nirenberg inequalities. In Section
5, we prove the weighted Rellich-type inequality and
Rellich-Sobolev inequality. In Section 6, we study the Hardy-type
inequalities with remainder term.
\section {Carnot group}

 A Carnot group (see \cite {3}, \cite{4},
\cite{18}, \cite{20}, \cite{21}, \cite {35} and \cite {40}) is a
connected, simply connected, nilpotent Lie group $\mathbb{G}$
whose Lie algebra $\mathcal{G}$ admits a stratification. That is,
there exist linear subspaces $V_1, ..., V_k$  of $\mathcal{G}$
such that
\begin{equation} \mathcal{G}=V_1\oplus...\oplus V_k, \quad [V_1, V_i]=V_{i+1},
\quad \text{for}\quad i=1,2, ..., k-1 \quad\text {and}\quad [V_1,
V_k]=0
\end{equation}
 where $[V_1, V_i]$ is the subspace  of $\mathcal{G}$
generated by the elements $[X,Y]$ with $X\in V_1$ and $Y\in V_i$.
This defines a $k$-step Carnot group and integer $k\ge 1$  is
called the step of $ \mathbb{G}$.

Via the exponential map, it is possible to induce on $\mathbb{G}$
a family of automorphisms of the group, called dilations,
$\delta_{\lambda}: \mathbb{R}^n\longrightarrow \mathbb{R}^n
(\lambda>0)$ such that
\[\delta_{\lambda}(x_1, ...,x_n)=(\lambda^{\alpha_1} x_1,...,\lambda^{\alpha_n}x_n)\]
where $1=\alpha_1=...=\alpha_m<\alpha_{m+1}\le ...\le \alpha_n$
are integers and $m=\text{dim}(V_1)$.

The group law can be written in the following form
\begin{equation}
x\cdot y=x+y+P(x,y), \quad x, y\in \mathbb{R}^n
\end{equation}
where $ P:\mathbb{R}^n\times \mathbb{R}^n\longrightarrow
\mathbb{R}^n$ has polynomial components and $P_1=...=P_m=0$. Note
that the inverse $x^{-1}$ of an element $x\in \mathbb{G}$ has the
form $x^{-1}=-x=(-x_1,..., -x_n)$.

Let $X_1, . . ., X_m$ be a family of left invariant vector fields
which form an orthonormal basis of $ V_1\equiv\mathbb{R}^m$ at the
origin, that is, $X_1(0)=\partial_{x_1}, . . ., X_m(0)=\partial
_{x_m}$. The vector fields $X_j$  have polynomial coefficients and
can be assumed to be of the form
\[X_j(x)=\partial_j+\sum_{i=j+1}^n a_{ij}(x)\partial_i, \quad
X_j(0)=\partial_j, j=1, . . ., m,\] where  each polynomial
$a_{ij}$ is homogeneous with respect to the dilations of the
group, that is
$a_{ij}(\delta_{\lambda}(x))=\lambda^{\alpha_i-\alpha_j}a_{ij}(x)$.
The horizontal  gradient on Carnot group $\mathbb{G}$ is the
vector valued operator
\[\nabla_{\mathbb{G}}=(X_1, . . ., X_m)\] where $X_1, . . ., X_m$ are the generators of $ \mathbb{G}$.
The sub-Laplacian is the second-order partial differential
operator on $\mathbb{G}$  given by
\[\Delta_{ \mathbb{G}}=\sum_{j=1}^m X_j^2.\]
The fundamental solution $u$ for $\Delta_{\mathbb{G}}$ is defined
to be a weak solution to the equation
\[-\Delta_{\mathbb{G}}u=\delta\] where $\delta$ denotes the Dirac distribution with singularity at
the neutral element $0$ of $ \mathbb{G}$. In \cite{18} Folland
proved that in any Carnot group $ \mathbb{G}$,  there exists a
homogeneous norm $N$ such that
\[u=N^{2-Q}\] is a fundamental solution for $\Delta_{
\mathbb{G}}$ ( see also \cite{7}).

We now set
\begin{equation}
N(x):=
\begin{cases}
u^{\frac{1}{2-Q}}&\quad\text{if}\quad x\neq 0,\\
 0 &\quad\text{if}\quad x=0.
 \end{cases}
 \end{equation}
 We recall
that a homogeneous norm on $ \mathbb{G}$ is a continuous
 function $N:\mathbb{G}\longrightarrow [0, \infty )$ smooth away from the origin
which satisfies the conditions : $N(\delta_{\lambda}(x))=\lambda
N(x)$, $N(x^{-1})=N(x)$ and $N(x)=0$ iff $x=0$.

 The curve $\gamma:[a,b]\subset \mathbb{R}\longrightarrow
\mathbb{G}$ is called horizontal if its tangents lie in $V_1$,
i.e,  $\gamma'(t)\in \text{\textit{span}}\{X_1, . . ., X_m\}$ for
all $t$. Then, the Carnot-Car\'ethedory distance $d_{CC}(x,y)$
between two points $x,y\in \mathbb{G}$ is defined to be the
infimum of all horizontal lengths $\int_a^b \langle \gamma'(t),
\gamma'(t)\rangle^{1/2} dt$ over all horizontal curves
$\gamma:[a,b]\longrightarrow \mathbb{G}$ such that $\gamma(a)=x$
and $\gamma(b)=y$. Notice that $d_{CC}$ is a homogeneous norm and
satisfies the invariance property

\[d_{CC}(z\cdot x, z\cdot y)=d_{CC}(x,y), \quad \text{for all} \, x,y,
z\in \mathbb{G},\] and is homogeneous of degree one with respect
to the dilation $\delta_{\lambda}$, i.e.
\[d_{CC}(\delta_{\lambda}(x), \delta_{\lambda}(y))=\lambda
d_{CC}(x,y), \quad \forall\, x, y, z \in \mathbb{G}, \text{for
all} \, \lambda>0.\]

The Carnot-Careth\'edory balls are defined by $B(x, R)=\{y\in
\mathbb{G} |d_{CC}(x,y)<R\}$. By left-translation and dilation, it
is easy to see that the Haar measure of $B(x, R)$ is proportional
by $R^Q$. More precisely
\[ |B(x, R)|=R^Q |B(x,1)|=R^Q |B(0,1)|\] where \[Q=\sum_{j=1}^k j (\text{dim}V_j)\] is the
homogeneous dimension of $ \mathbb{G}$.

It is well known that  Sobolev inequalities are important in the
study of partial differential equations, especially in the study
of those arising from geometry and physics. The following Sobolev
inequality holds  on $ \mathbb{G}$ \cite {18}
\begin{equation}
\Big(\int_{\mathbb{G}}|\nabla_{\mathbb{G}}\phi(x)|^2dx\Big)^{1/2}\ge
C
\Big(\int_{\mathbb{G}}|\phi(x)|^{\frac{2Q}{Q-2}}dx\Big)^{\frac{Q-2
}{2Q}}
\end{equation}
( see also for weighted and higher order extensions \cite
{32,33,34} ). It is a more difficult problem to determine the
sharp constant $C$ in (2.4)  on general Carnot groups. The only
results that have so far been proven are in the case of Heisenberg
group $\mathbb{H}^n$ by Jerison and Lee \cite{29} and Iwasawa-type
groups ( a particular sub-class of $\textit{H}$-type groups) by
Garofalo and Vassilev \cite {24}. We should mention that the sharp
constants in \cite{29} and \cite{24} lead us to obtain explicit
constant in Corollary 4.3.

\section{Hardy-type inequalities on Carnot groups of step 2}
Among  Carnot groups of step two,  the Heisenberg group  and
Heisenberg type (\textit{H}-type) groups are of particular
significance. These groups appear naturally in analysis, geometry,
representation theory and mathematical physics. In this section,
we first prove Hardy-type inequalities on the Heisenberg group and
we extend this result to the \textit{H}-type group.

\noindent {\textbf{Heisenberg group}.} The Heisenberg group $
\mathbb{H}^n$ is an example of a noncommutative Carnot group.
Denoting points in $ \mathbb{H}^n$ by $(z,t)$ with $z=(z_1, . . .,
z_n)\in \mathbb{C}^n$ and $t\in \mathbb{R}$ we have the group law
given as
\[(z,t)\circ (z',t')=(z+z', t+t'+2\sum_{j=1}^n
Im(z_j\bar{z}_j'))\] With the notation $z_j=x_j+iy_j$, the
horizontal space $V_1$ is spanned by the basis
\[X_j=\frac{\partial}{\partial x_j}
+2y_j\frac{\partial }{\partial t}\quad\text{and} \quad
Y_j=\frac{\partial}{\partial y_j}-2x_j\frac{\partial }{\partial
t}.\] The one dimensional center $V_2$ is spanned by the vector
field $T=\frac{\partial}{\partial t}$.  We have the commutator
relations $[X_j, Y_j]=-4T$, and all other brackets of $\{X_1, Y_1,
. . ., X_n, Y_n\}$ are zero. The sub-elliptic gradient is the $2n$
dimensional vector field given by
\[\nabla_{ \mathbb{H}^n}=(X_1, . . ., X_n, Y_1, . . ., Y_n)\] and
the Kohn Laplacian on $ \mathbb{H}^n$ is the operator \[\Delta_{
\mathbb{H}^n}=\sum_{j=1}^n (X_j^2+Y_j^2).\] A homogeneous norm on
$ \mathbb{H}^n$ is given by
\[\rho=|(z,t)|=(|z|^4+t^2)^{1/4}\] and the homogeneous dimension of $\mathbb{H}^n$ is $Q=2n+2$.

A remarkable analogy between Kohn Laplacian and the classical
Laplace operator has been obtained by Folland \cite{17}. He found
that the fundamental solution of $-\Delta_{ \mathbb{H}^n}$ with
pole at zero is given by
\[\Psi(z,t)=\frac{c_Q}{\rho(z,t)^{Q-2}}\quad\text{where}\quad
c_Q=\frac{2^{(Q-2)/2}\Gamma((Q-2)/4)^2}{\pi^{Q/2}}.\]

We now prove the following theorem  on the Heisenberg group $
\mathbb{H}^n$. In the various integral inequalities below (Section
3 and Section 4), we allow the values of the integrals on the
left-hand sides to be $+\infty$.  Before we proceed, we should
emphasize that the constant $C(Q,\alpha)=(\frac{Q+\alpha-2}{2})^2$
obtained in Section 3 and Section 4 is  sharp in the sense that if
it is replaced by an grater number the inequality fails.

\begin{theorem}
Let $\alpha\in \mathbb{R}$ and $\phi\in
C_0^{\infty}(\mathbb{H}^n\setminus \{0\})$. Then we have :
\[\int_{
\mathbb{H}^n}\rho^{\alpha}|\nabla_{\mathbb{H}^n}\phi|^2dzdt\ge
\Big(\frac{Q+\alpha-2}{2}\Big)^2\int_
{\mathbb{H}^n}\rho^{\alpha}\frac{|z|^2}{\rho^4}\phi^2dzdt\] where
$\rho=(|z|^4+l^2)^{1/4}$ is the homogeneous norm on $
\mathbb{H}^n$. Moreover, the constant $(\frac{Q+\alpha-2}{2})^2$
is sharp provided $Q+\alpha-2>0$.
\end{theorem}
\proof Let $\phi=\rho^{\beta}\psi$ where $\beta\in
\mathbb{R}\setminus \{0\}$ and $\psi \in
C_0^{\infty}(\mathbb{H}^n\setminus \{0\})$. A direct calculation
shows that
\begin{equation}
\rho^{\alpha}|\nabla_{ \mathbb{H}^n}\phi|^2
=\beta^2\rho^{\alpha+2\beta-2}|\nabla_{\mathbb{H}^n}\rho|^2\psi^2+
2\beta\rho^{\alpha+2\beta-1}\psi\nabla_{
\mathbb{H}^n}\rho\cdot\nabla_{
\mathbb{H}^n}\psi+\rho^{\alpha+2\beta}|\nabla_{\mathbb{H}}\psi|^2.
\end{equation}
It is easy to see that
\[|\nabla_{ \mathbb{H}^n}\rho|^2=\frac{|z|^2}{\rho^2}\]
and integrating (3.1) over $ \mathbb{H}^n$, we get
\begin{equation}
\begin{aligned}\int_{ \mathbb{H}^n}\rho^{\alpha}|\nabla_{
\mathbb{H}^n}\phi|^2dzdt&=\int_{
\mathbb{H}^n}\beta^2\rho^{\alpha+2\beta-4}|z|^2\psi^2dzdt+\int_{
\mathbb{H}^n} 2\beta\rho^{\alpha+2\beta-1}\psi\nabla_{
\mathbb{H}^n}\rho\cdot\nabla_{ \mathbb{H}^n}\psi dzdt\\&+\int_{
\mathbb{H}^n}\rho^{\alpha+2\beta}|\nabla_{\mathbb{H}^n}\psi|^2dzdt\\
\end{aligned}
\end{equation}
Applying integration by parts to the middle integral on the
right-hand side of (3.2), we obtain

\begin{equation}
\begin{aligned}
\int_{ \mathbb{H}^n}\rho^{\alpha}|\nabla_{
\mathbb{H}^n}\phi|^2dzdt&=\int_{
\mathbb{H}^n}\beta^2\rho^{\alpha+2\beta-4}|z|^2\psi^2dzdt-\frac{\beta}{\alpha+2\beta}
\int_{ \mathbb{H}^n}\Delta_{
\mathbb{H}^n}(\rho^{\alpha+2\beta})\psi^2dzdt\\&+\int_{
\mathbb{H}^n}\rho^{\alpha+2\beta}|\nabla_{\mathbb{H}^n}\psi|^2dzdt.
\end{aligned}
\end{equation}
One can show that
\begin{equation}\Delta_{
\mathbb{H}^n}(\rho^{\alpha+2\beta})=|z|^2\rho^{\alpha+2\beta-4}(\alpha+2\beta)(\alpha+2\beta+Q-2).
\end{equation}
Substituting (3.4) into (3.3) gives the following

\[\begin{aligned} \int_{
\mathbb{H}^n}\rho^{\alpha}|\nabla_{ \mathbb{H}^n}\phi|^2dzdt &=
(\beta^2-\beta(\alpha+2\beta+Q-2))\int_{
\mathbb{H}^n}\rho^{\alpha+2\beta-4}|z|^2\psi^2dzdt+ \int_{
\mathbb{H}^n}
\rho^{\alpha+2\beta}|\nabla_{\mathbb{H}^n}\psi|^2dzdt\\
&\ge(-\beta^2-\beta(\alpha+Q-2))\int_{
\mathbb{H}^n}\rho^{\alpha+2\beta-4}|z|^2\psi^2dzdt.
\end{aligned}
\]
Note that the function $\beta\longrightarrow
-\beta^2-\beta(\alpha+Q-2)$ attains the maximum for
$\beta=\frac{2-\alpha-Q}{2}$, and this maximum is equal to
$(\frac{Q+\alpha-2}{2})^2$. Therefore we have the following
inequality
\begin{equation}
\int_{ \mathbb{H}^n}\rho^{\alpha}|\nabla_{
\mathbb{H}^n}\phi|^2dzdt \ge
\Big(\frac{Q+\alpha-2}{2}\Big)^2\int_{
\mathbb{H}^n}\rho^{\alpha}\frac{|z|^2}{\rho^4}\phi^2dzdt.
\end{equation}

It only remains to show that the constant
$(\frac{Q+\alpha-2}{2})^2$ is sharp. The method  employed here is
quite standard which is adapted from the Euclidean case (see \cite
{22}). We now give proof for the Heisenberg group case and proof
for the H-type groups is similar. Let $\phi_{\epsilon}(z,t)$ be
the family of functions defined by
\begin{equation}
\phi_{\epsilon}(z,t)=
\begin{cases}
 1 &\quad\text{if}  \quad \rho \in [0,1],\\
\rho^{-(\frac{Q+\alpha-2}{2}+\epsilon)} &\quad \text{if} \quad
\rho
>1,
\end{cases}
\end{equation}
where $\epsilon>0$ and $\rho=|(z,t)|=(|z|^4+t^2)^{1/4}$. It
follows that
\[|\rho^{\alpha}\nabla_{\mathbb{H}^n}\phi_{\epsilon}(z,t)|^2=\Big(\frac{Q+\alpha-2}{2}+\epsilon\Big)^2|z|^2\rho^{-(Q+2+2\epsilon)}.\]
In the sequel we indicate  $B_1=\{(z,t):\rho\le 1\}$ $\rho$-ball
centered at the origin in $ \mathbb{H}^n$ with radius $1$.

 By direct computation we get \begin{equation}\begin{aligned}
\int_{\mathbb{H}^n}
\rho^{\alpha}\frac{|z|^2}{\rho^4}\phi_{\epsilon}^2dzdt&=\int_{B_1}\rho^{\alpha-4}|z|^2\phi_{\epsilon}^2dzdt
+\int_{\mathbb{H}^n\setminus
B_1}\rho^{\alpha-4}|z|^2\phi_{\epsilon}^2 dzdt\\
&=\int_{B_1}\rho^{\alpha-4}|z|^2dzdt+\int_{\mathbb{H}^n\setminus
B_1}|z|^2\rho^{-(Q+2+2\epsilon)}dzdt\\
&=\int_{B_1}\rho^{\alpha-4}|z|^2dzdt+
(\frac{Q+\alpha-2}{2}+\epsilon)^{-2}\int_{\mathbb{H}^n}\rho^{\alpha}
|\nabla_{\mathbb{H}^n}\phi_{\epsilon}|^2dzdt.\end{aligned}\end{equation}
Since $Q+\alpha-2>0$ then the first integral on the right hand
side of (3.7) is integrable and we conclude by
$\epsilon\longrightarrow 0$.

\endproof
\medskip
The following theorem shows that the weight function
$\rho^{\alpha}$ has a significant effect on the sharp constant
$(\frac{Q+\alpha-2}{2})^2$ whereas the new weight function
$|\nabla_{\mathbb{H}^n}\rho|^{\gamma}$ has no effect.
\begin{theorem}
Let $\alpha, \gamma \in \mathbb{R}$ and $\phi\in
C_0^{\infty}(\mathbb{H}^n\setminus \{0\})$. Then we have :
\[\int_{
\mathbb{H}^n}\rho^{\alpha}|\nabla_{\mathbb{H}^n}\rho|^{\gamma}
|\nabla_{\mathbb{H}^n}\phi|^2dzdt\ge (\frac{Q+\alpha-2}{2})^2\int_
{\mathbb{H}^n}\rho^{\alpha-2}
|\nabla_{\mathbb{H}^n}\rho|^{\gamma+2}\phi^2dzdt\] where
$\rho=(|z|^4+l^2)^{1/4}$ is the homogeneous norm on $
\mathbb{H}^n$. Moreover, the constant $(\frac{Q+\alpha-2}{2})^2$
is sharp provided $Q+\alpha-2>0$.
\end{theorem}
\proof The proof is similar to the proof of Theorem 3.1. We only
need to note that \[\nabla_{\mathbb{H}^n}\rho\cdot
\nabla_{\mathbb{H}^n}(|\nabla_{\mathbb{H}^n}\rho|^{\gamma})=0.
\]
\endproof

\noindent{\textbf{Heisenberg-type group}.} Another important model
of Carnot groups are the \textit{H}-type (Heisenberg type) groups
which were introduced by Kaplan \cite {30} as direct
generalizations of the Heisenberg group $\mathbb{H}^n$. An
\textit{H}-type group is a Carnot group with a two-step Lie
algebra $ \mathcal{G}=V_1\oplus V_2$ and an inner product
$\langle, \rangle$ in $ \mathcal{G}$ such that the linear map
\[ J:V_2\longrightarrow \text{End}V_1,\] defined by the condition
\[\langle J_z(u), v\rangle=\langle z, [u,v]\rangle,\quad u,v\in V_1, z\in V_2\] satisfies
\[J_z^2=-||z||^2\mathbf{Id}\]
for all $z\in V_2$, where $||z||^2= \langle z,z\rangle$.

Sub-Laplacian is defined in terms of a fixed basis $X_1,. . . ,
X_m$ for $V_1$:

\begin{equation}
\Delta_\mathbb{G}=\sum_{i=1}^mX_i^2.
\end{equation}

The exponential mapping of a simply connected Lie group is an
analytic diffeomorphism. One can then define analytic mappings
$v:\mathbb{G} \longrightarrow V_1$ and
$z:\mathbb{G}\longrightarrow V_2$ by
\[x=\text{exp}(v(x)+z(x))\] for every $x\in \mathbb{G}$.
In \cite {30} Kaplan proved that there exists a constant $c>0$
such that the function
\[\Phi (x)=c\Big(|v(x)|^4+16|z(x)|^2\Big)^{\frac{2-Q}{4}}\]
is a  fundamental solution for the operator $ -\Delta_\mathbb{G}$.
We note that
\begin{equation}
K(x)=\Big(|v(x)|^4+16|z(x)|^2\Big)^{\frac{1}{4}}
\end{equation}
defines a homogeneous norm  and $Q=m+2k$ is the homogeneous
dimension of $\mathbb{G}$ where $m=$dim$V_1$ and $k=$dim$V_2$.
This result generalized Folland's fundamental solution for the
Heisenberg group $ \mathbb{H}^n$ \cite{17}. It is remarkable that
the homogeneous norm $K(x)$  is involved also in the expression of
the fundamental solution of the following $p-$sub-Laplace operator
\begin{equation}
\mathcal{L}_pu=\sum_{i=1}^mX_i(|Xu|^{p-2}X_iu), \quad 1<p<\infty.
\end{equation}
More precisely, Capogna, Danielli and Garofalo \cite{11} proved
that for every $1<p<\infty$ there exists $c_p>0$  such that the
function
\begin{equation}
\Gamma_p(x)=
\begin{cases}
c_pK^{(p-Q)/(p-1)} &\quad \text{when}\quad p\ne Q,\\
-c_p\log K&\quad\text{when}\quad p=Q,
\end{cases}
\end{equation}
is a fundamental solution for the operator $-\mathcal{L}_p$ (see
also \cite {26} for the case $p=Q$).

We cite, without proof of the following, useful formulas which can
be found in \cite{11}. Let $u$ be a radial function, i.e.,
$u(x)=f(K(x))$ where $f\in C(\mathbb{R})$ then

\[ |\nabla_{ \mathbb{G}}u|^2=\frac{|v|^2}{K^2}|f'(K)|^2.\] Moreover  if  $u(x)=f(K(x))$  and $f\in
C^2(\mathbb{R})$ then
\begin{equation}
\begin{aligned}
\Delta_{\mathbb{G}}u&=|\nabla_{\mathbb{G}}K(x)|^2\Big[f''(K)+\frac{Q-1}{K}f'(K)\Big]\\
&=\frac{v^2}{K^2} \Big[f''(K)+\frac{Q-1}{K}f'(K)\Big]
\end{aligned}
\end{equation}
 at every point $x\in \mathbb{G}\setminus \{0\}$ where
$f'(K(x))\neq 0$.

Another important fact that $K(x)$ satisfies the so-called
$\infty$-sub-Laplace  equation :
\[ \Delta_{ \mathbb{G}, \infty}K=\frac{1}{2}\langle
\nabla_{\mathbb{G}}(|\nabla_{\mathbb{G}}K|^2),
\nabla_{\mathbb{G}}K \rangle=0\] at every point $x\in
\mathbb{G}\setminus \{0\}$. (See \cite{30} and \cite{11} for
further information on \textit{H}-type groups)
\medskip

We now have the following theorem on the \textit{H}-type group :

\begin{theorem}
Let $ \mathbb{G}$ be an \textit{H}-type group with homogeneous
dimension $Q=m+2k$ and let $\alpha, \gamma \in \mathbb{R}$ and
$\phi\in C_0^{\infty}(\mathbb{G}\setminus \{0\})$. Then the
following inequality is valid :
\begin{equation}\int_{
\mathbb{G}}K^{\alpha}|\nabla_{\mathbb{G}}K|^{\gamma}|\nabla_{\mathbb{G}}\phi|^2dx\ge
(\frac{Q+\alpha-2}{2})^2\int_
{\mathbb{G}}K^{\alpha-2}|\nabla_{\mathbb{G}}K|^{\gamma+2}
\phi^2dx\end{equation} where $K(x)=(|v(x)|^4+16|z(x)|^2)^{1/4}$.
Moreover, the constant $(\frac{Q+\alpha-2}{2})^2$ is sharp
provided $Q+\alpha-2>0$.
\end{theorem}
\proof The proof is identical to the Heisenberg group case.
\endproof

\section{Hardy-type inequalities on Carnot groups of arbitrary step }
In this section,  we consider Carnot group $ \mathbb{G}$ of any
step $k$ with the homogeneous norm $N=u^{1/(2-Q)}$ associated to
Folland's solution $u$ for the sub-Laplacian $\Delta_{
\mathbb{G}}$ \cite{18}. We have the following theorem:
\begin{theorem} Let $\mathbb{G}$ be a Carnot group with homogeneous dimension $Q\ge 3$  and  let $\phi\in
C_0^{\infty}(\mathbb{G}\setminus \{0\})$, $\alpha\in \mathbb{R}$,
$Q+\alpha-2>0$. Then the following inequality is valid
\begin{equation}\int_{\mathbb{G}} N^{\alpha}|\nabla_{\mathbb{G}} \phi|^2dx
\ge \Big(\frac{Q+\alpha-2}{2}\Big)^2 \int_{\mathbb{G}}
N^{\alpha}\frac{|\nabla_{\mathbb{G}} N|^2}{N^2}\phi
^2dx.\end{equation}  Furthermore,  the constant
$C(Q,\alpha)=(\frac{Q+\alpha-2}{2})^2$ is sharp.
\end{theorem}
\proof

 Let $\phi=N^{\beta}\psi$  where $\psi\in
C_0^{\infty}(\mathbb{G}\setminus \{0\})$ and $\beta\in
\mathbb{R}\setminus \{0\}$. A direct calculation
 shows that
 \begin{equation}
|\nabla_{\mathbb{G}}(N^{\beta}\psi)|^2=\beta^2N^{2\beta-2}|\nabla_{\mathbb{G}}
N|^2\psi^2+2\beta N^{2\beta-1}\psi\nabla_{\mathbb{G}} N\cdot
\nabla_{\mathbb{G}} \psi+N^{2\beta}|\nabla_{\mathbb{G}}\psi|^2.
\end{equation}
Multiplying both sides of (4.2) by the $N^{\alpha}$ and applying
integration by parts over $ \mathbb{G}$ gives

\begin{equation}\begin{aligned}
\int_{\mathbb{G}}N^{\alpha}|\nabla_{\mathbb{G}}\phi|^2dx
&=\beta^2\int_{\mathbb{G}} N^{\alpha+2\beta-2}|\nabla_{\mathbb{G}}
N|^2\psi ^2dx-\frac{\beta}{\alpha+2\beta}\int_{\mathbb{G}}
\Delta_{\mathbb{G}}(N^{\alpha+2\beta})\psi^2dx\\&+\int_{\mathbb{G}}N^{\alpha+2\beta}|\nabla_{\mathbb{G}}\psi|^2dx\\
& \ge \beta^2\int_{\mathbb{G}}
N^{\alpha+2\beta-2}|\nabla_{\mathbb{G}} N|^2\psi
^2dx-\frac{\beta}{\alpha+2\beta}\int_{\mathbb{G}}
\Delta_{\mathbb{G}}(N^{\alpha+2\beta})\psi^2dx.
\end{aligned}
\end{equation}
A straightforward calculation shows that

\begin{equation}
-\frac{\beta}{\alpha+2\beta} \Delta_{
\mathbb{G}}(N^{\alpha+2\beta})=-\beta(\alpha+2\beta+Q-2)N^{\alpha+2\beta-2}|\nabla_{\mathbb{G}}N|^2
-\frac{\beta}{2-Q}N^{\alpha+2\beta+Q-2}\Delta_{\mathbb{G}}u.
\end{equation}
Substituting (4.4) into (4.3) and  using the fact that
$\psi^2=N^{-2\beta}\phi^2$, we get the following :

\[
\int_{\mathbb{G}}N^{\alpha}|\nabla_{\mathbb{G}}\phi|^2dx\ge (
-\beta^2-\beta(\alpha+Q-2)\int_{
\mathbb{G}}N^{\alpha}\frac{|\nabla_{ \mathbb{G}}
N|^2}{N^2}\phi^2dx-\frac{\beta}{2-Q}\int_{
\mathbb{G}}(\Delta_{\mathbb{G}}u)N^{\alpha+Q-2}\phi^2dx.
\]
Since $u$ is the fundamental solution of sub-Laplacian
$\Delta_{\mathbb{G}}$ on Carnot group $\mathbb{G}$, we get
\[-\int_{
\mathbb{G}}(\Delta_{\mathbb{G}}u)N^{\alpha+Q-2}\phi^2dx=
N^{\alpha+Q-2}(0)\phi^2(0)=0.\]  We now obtain

\[\int_{\mathbb{G}}N^{\alpha}|\nabla_{\mathbb{G}}\phi|^2dx\ge (
-\beta^2-\beta(\alpha+Q-2)\int_{
\mathbb{G}}N^{\alpha}\frac{|\nabla_{ \mathbb{G}}
N|^2}{N^2}\phi^2dx.\] Choosing
\[\beta=\frac{2-Q-\alpha}{2}\]  gives the following
inequality
\begin{equation}
\int_{\mathbb{G}} N^{\alpha}|\nabla_{\mathbb{G}} \phi|^2dx \ge
\Big(\frac{Q+\alpha-2}{2}\Big)^2 \int_{\mathbb{G}}
N^{\alpha}\frac{|\nabla_{\mathbb{G}} N|^2}{N^2}\phi ^2dx.
\end{equation}

To show that the constant $\Big(\frac{Q+\alpha-2}{2}\Big)^2$ is
sharp, we use the following family of functions
\begin{equation} \phi_{\epsilon}(x)=
\begin{cases}
 1 &\quad\text{if}  \quad N(x)\in [0,1],\\
N^{-(\frac{Q+\alpha-2}{2}+\epsilon)} &\quad \text{if} \quad
N(x)>1,
\end{cases}
\end{equation}
and pass to the limit as $\epsilon\longrightarrow 0$. We should
indicate that same test function lead us to obtain sharp constant
in Theorem (4.3).  Here we notice that $|\nabla_{\mathbb{G}}N|$ is
uniformly bounded and polar coordinate integration formula holds
 on  $\mathbb{G}$ (\cite{20}).
\endproof

\noindent\textbf{Remark 4.1.} In the abelian case, when
$\mathbb{G}=\mathbb{R}^n$ with the ordinary dilations, one has
$\mathcal{G}=V_1=\mathbb{R}^n$ so that $Q=n$. Now it is clear that
the inequality (4.1)  with the homogeneous norm $N(x)=|x|$ and
$\alpha=0$ recovers the Hardy inequality (1.1).
\medskip

\noindent\textbf{Uncertainty Principle Inequality.} The classical
uncertainty principle was developed in the context of quantum
mechanics by Heisenberg \cite {28}. It says that the position and
momentum of a particle cannot be determined exactly at the same
time but only with an ``uncertainty". The harmonic analysis
version of uncertainty principle states that a function on the
real line and its Fourier transform can not be simultaneously well
localized. It has been widely studied in quantum mechanics and
signal analysis. There are various forms of the uncertainty
principle. For an overview we refer to Folland's and Sitaram's
paper \cite {19}.

The uncertainty principle on the Euclidean space $ \mathbb{R}^n$
can be  stated in the following way:
\begin{equation}
\Big(\int_{\mathbb{R}^n} |x|^2
|f(x)|^2dx\Big)\Big(\int_{\mathbb{R}^n} |\nabla f(x)|^2 dx\Big)\ge
\frac{n^2}{4} \Big(\int_{\mathbb{R}^n} |f(x)|^2 dx \Big)^2
\end{equation}
for all $f\in L^2( \mathbb{R}^n)$. An analogue of the above
inequality (4.7) for the Heisenberg group $ \mathbb{H}^n$ was
established by Garofalo and Lanconelli \cite {23}. Thangavelu
\cite {41}, and Sitaram, Sundari and Thangavelu \cite{39} have
also obtained related but inequivalent analogues of Heisenberg's
inequality for functions on the Heisenberg group $ \mathbb{H}^n$.

In the following corollaries,  we present the analogues of (4.7)
for general Carnot groups. The proof of the corollaries is based
on the Hardy-type inequality (4.1) and the Cauchy-Schwarz
inequality.  We should mention that  Ciatti, Ricci and Sundari
\cite {12} have also obtained a version of uncertainty principle
inequality on nilpotent stratified Lie groups of step two (Carnot
group of step 2) which is not equivalent to our result (see also
Corollary 5.1 ).

\begin{corollary} Let $\mathbb{G}$ be a
Carnot group with homogeneous dimension $Q\ge 3$. Then for every
$\phi\in C_0^{\infty}(\mathbb{G}\setminus \{0\})$
\begin{equation}
\Big(\int_{\mathbb{G}}\frac{
N^2}{|\nabla_{\mathbb{G}}N|^2}\phi^2dx\Big)\Big(\int_{\mathbb{G}}|\nabla_{\mathbb{G}}
\phi|^2dx\Big)\ge
\Big(\frac{Q-2}{2}\Big)^2\Big(\int_{\mathbb{G}}\phi^2 dx\Big)^2.
\end{equation}
\end{corollary}

\begin{corollary} Let $\mathbb{G}$ be a
Carnot group with homogeneous dimension $Q\ge 3$. Then for every
$\phi\in C_0^{\infty}(\mathbb{G}\setminus \{0\})$
\begin{equation}
\Big(\int_{\mathbb{G}}
N^2|\nabla_{\mathbb{G}}N|^2\phi^2dx\Big)\Big(\int_{\mathbb{G}}|\nabla_{\mathbb{G}}
\phi|^2dx\Big)\ge
\Big(\frac{Q-2}{2}\Big)^2\Big(\int_{\mathbb{G}}|\nabla_{\mathbb{G}}N|^2\phi^2
dx\Big)^2.
\end{equation}
\end{corollary}
\medskip

\noindent\textbf{Remark 4.2.} It is well known that equality is
attained in  uncertainty principle inequality (4.7) only for
Gaussian functions. As mentioned above, this fact has been also
extended to the Heisenberg group by  Garofalo and Lanconelli
\cite{ 23}. It is natural to search an analogue of this phenomena
for general Carnot groups. We should notice that with
$\frac{Q-2}{2}$ replaced by $\frac{Q}{2}$ equality is attained in
Corollary (4.8) and Corollary (4.9) if $\phi(x)=C e^{-\beta
N^2(x)}$ for some $C\in \mathbb{R}, \beta>0$. ( Note that
$\nabla_{\mathbb{G}}N(x)\neq 0$ for (Haar) a.e $x\in \mathbb{G}$
and $|\nabla_{\mathbb{G}}N|$  is uniformly bounded on $
\mathbb{G}$ \cite {3}.)
\medskip

The following corollaries are the  consequence of the Hard-type
(4.1) and Sobolev (2.4) inequalities. These inequalities are
extensions of the Caffarelli-Kohn-Nirenberg \cite{10} inequality
to Carnot groups.

\begin{corollary} Let $\mathbb{G}$ be a
Carnot group with homogeneous dimension $Q\ge 3$  and $0\le s \le
2$. Then for every $\phi\in C_0^{\infty}(\mathbb{G}\setminus
\{0\})$ then there exists a constant $C>0$ such that
\[
\int_{\mathbb{G}}|\nabla_{\mathbb{G}}\phi|^2dx\ge C\Big(
\int_{\mathbb{G}}\Big(\frac{|\nabla_{ \mathbb{G}}N|}{N}\Big)^s
|\phi|^{2(\frac{Q-s}{Q-2})}dx\Big)^{\frac{Q-2}{Q-s}}.
\]
\end{corollary}
\medskip

We now have the following weighted inequality on metric ball $B$.
The proof of this is inequality based on the Hardy-type and
weighted Sobolev inequalities \cite {34}. We note that the weight
function $N^{\alpha}$ in Corollary 4.4 satisfies the Muckenhoupt
$A_2$ condition and other requirements for the existence of
weighted Sobolev inequality. We recall that a weight $w(x)$
satisfies Muckenhoupt $A_p$ condition for $1<p<\infty$ if there is
a constant $C$ such that
\[\Big(\frac{1}{|B|}\int_{B}w(x)dx\Big)^{1/p}\Big(\frac{1}{|B|}\int_{B}w(x)^{-p'/p}dx\Big)^{\frac{1}{p'}}\le
C\] for all metric balls $B$.  If $w(x)\in A_p$ then we have  $
w(x)^{-p'/p}\in A_{p'}$ where $p'$ is the dual exponent to $p$
given by $\frac{1}{p}+\frac{1}{p'}=1$.
\begin{corollary} Let $\mathbb{G}$ be a
Carnot group with homogeneous dimension $Q\ge 3$  and let $B$ be a
metric ball in  $\mathbb{G}$, $2-Q<\alpha<Q$, $0\le s \le 2$. Then
for every $\phi\in C_0^{\infty}(B\setminus \{0\})$ then there
exists a constant $C>0$ such that
\[
\int_{B}N^{\alpha}|\nabla_{\mathbb{G}}\phi|^2dx\ge C\Big(
\int_{B}N^{\alpha}\Big(\frac{|\nabla_{ \mathbb{G}}N|}{N}\Big)^s
|\phi|^{2(\frac{Q-s}{Q-2})}dx\Big)^{\frac{Q-2}{Q-s}}.
\]
\end{corollary}

\medskip
\noindent \textbf{Polarizable Carnot group.} The second main
result of this section is to establish weighted Hardy-type
inequality including the weight function  $\
|\nabla_{\mathbb{G}}N|^{\gamma}$  as in the Section 3. We should
indicate that Hardy-type inequality in \cite {14} on Carnot groups
does not include the weight function $
|\nabla_{\mathbb{G}}N|^{\gamma}$. We now establish such a
inequality on polarizable Carnot groups. This class of groups were
introduced by Balogh and Tyson \cite{3} and admit the analogue of
polar coordinates.

A Carnot group $\mathbb{G}$  is said to be polarizable if the
homogeneous norm $N=u^{1/(2-Q)}$,
 associated to Folland's solution $u$ for the sub-Laplacian
 $\Delta_{\mathbb{G}}$,   satisfies the following $\infty$- sub-Laplace
 equation,

 \begin{equation}
\Delta_{\mathbb{G},\infty}N:= \frac{1}{2}\langle\nabla_{
\mathbb{G}}(|\nabla_{ \mathbb{G}}N|^2), \nabla_{
\mathbb{G}}N\rangle=0, \quad\quad \text{in} \quad
\mathbb{G}\setminus \{0\}.
 \end{equation}

It has been proved that the homogeneous norm (3.9) satisfies the
equation (4.10) ( see  \cite {3} and \cite {15}) . This
 result implies that H-type groups are polarizable Carnot
 groups. Unfortunately, at the present time,  it is unknown an example of polarizable
 Carnot group which is not of \textit{H}-type.

Balogh and Tyson \cite{3} proved  that the homogeneous norm
$N=u^{1/(2-Q)}$, associated to Folland's solution $u$ for the
sub-Laplacian
 $\Delta_{\mathbb{G}}$, enters also in the expression of the fundamental solution of the sub-elliptic $p$-Laplacian
\begin{equation}
\Delta_{ \mathbb{G}, p}u=\sum_{i=1}^mX_i(|Xu|^{p-2}X_iu), \quad
1<p<\infty.
\end{equation}

More precisely, Balogh and Tyson \cite{3} proved that for every
$1<p<\infty, p=Q$ there exists $c_p>0$ such that the fundamental
solution of $-\Delta_{\mathbb{G},p}$
 is given
by
\begin{equation}
u_p=
\begin{cases}
\begin{aligned}  c_pN^{\frac{p-Q}{p-1}}, &\quad\text{if}\quad p \neq
Q,\\
-c_Q\text{log}N, &\quad\text{if}\quad p=Q.
\end{aligned}
\end{cases}
\end{equation}

We are now ready to state our the second main theorem in this
section.

\begin{theorem} Let $\mathbb{G}$ be a polarizable Carnot group with homogeneous norm $N=u^{1/(2-Q)}$ and  let $\phi\in
C_0^{\infty}(\mathbb{G}\setminus \{0\})$, $\alpha \in
\mathbb{R}$,$\gamma>-1$,  $Q\ge 3$, $Q+\alpha-2>0$. Then the
following inequality is valid
\begin{equation}\int_{\mathbb{G}} N^{\alpha}|\nabla_{\mathbb{G}} N|^{\gamma}|\nabla_{\mathbb{G}} \phi|^2dx
\ge \Big(\frac{Q+\alpha-2}{2}\Big)^2 \int_{\mathbb{G}}
N^{\alpha}\frac{|\nabla_{\mathbb{G}} N|^{\gamma+2}}{N^2}\phi
^2dx.\end{equation} Furthermore, the constant
$C(Q,\alpha)=(\frac{Q+\alpha-2}{2})^2$ is  sharp in the sense that
if it is replaced by an grater number the inequality fails.
\end{theorem}
\proof

 Let $\phi=N^{\beta}\psi$  where $\psi\in
C_0^{\infty}(\mathbb{G}\setminus \{0\})$ and $\beta\in
\mathbb{R}\setminus \{0\}$. A direct calculation
 shows that
 \begin{equation}
 \begin{aligned}
\int _{\mathbb{G}} N^{\alpha}|\nabla_{\mathbb{G}}
N|^{\gamma}|\nabla_{\mathbb{G}} \phi|^2dx &= \beta^2 \int
_{\mathbb{G}}N^{\alpha+2\beta-2}|\nabla_{\mathbb{G}}
N|^{\gamma+2}\psi^2 dx\\ &+ 2\beta\int _{\mathbb{G}}
N^{\alpha+2\beta-1}|\nabla_{\mathbb{G}}
N|^{\gamma}\psi\nabla_{\mathbb{G}} N\cdot \nabla_{\mathbb{G}} \psi
dx \\ &+ \int _{\mathbb{G}}N^{\alpha+2\beta}|\nabla_{\mathbb{G}}
N|^{\gamma}|\nabla_{\mathbb{G}}\psi|^2dx. \end{aligned}
\end{equation}
Applying integration by parts to the middle term:\\
\begin{equation}
 \begin{aligned}
\int _{\mathbb{G}} N^{\alpha}|\nabla_{\mathbb{G}}
N|^{\gamma}|\nabla_{\mathbb{G}} \phi|^2dx &= \beta^2 \int
_{\mathbb{G}}N^{\alpha+2\beta-2}|\nabla_{\mathbb{G}}
N|^{\gamma+2}\psi^2 dx\\ & -\beta\int _{\mathbb{G}}\psi^2
\text{div}\big( N^{\alpha+2\beta-1}|\nabla_{\mathbb{G}}
N|^{\gamma}\nabla_{\mathbb{G}} N\big) dx
\\ &+ \int _{\mathbb{G}}N^{\alpha+2\beta}|\nabla_{\mathbb{G}}
N|^{\gamma}|\nabla_{\mathbb{G}}\psi|^2dx. \end{aligned}
\end{equation}
 We now choose $\gamma=p-2>1$ and
$ \alpha+2\beta-1=1-Q$, we get
\begin{equation}
\begin{aligned}
\int _{\mathbb{G}} N^{\alpha}|\nabla_{\mathbb{G}}
N|^{\gamma}|\nabla_{\mathbb{G}} \phi|^2dx &=
\Big(\frac{Q+\alpha-2}{2}\Big)^2  \int
_{\mathbb{G}}N^{\alpha-2}|\nabla_{\mathbb{G}}
N|^{\gamma+2}\phi^2 dx\\
&-\beta c_p^{1-p}\int _{\mathbb{G}}(\Delta_{\mathbb{G},
p}(u_p))N^{-2\beta}\phi^2 dx+ \int
_{\mathbb{G}}N^{2-Q}|\nabla_{\mathbb{G}}
N|^{\gamma}|\nabla_{\mathbb{G}}\psi|^2dx.
\end{aligned}
\end{equation}
Since $u_p$ is the fundamental solution of sub-p-Laplacian
$-\Delta_{\mathbb{G},p}$, we get
\[-\int _{\mathbb{G}}(\Delta_{\mathbb{G},
p}(u_p))N^{Q+\alpha-2}\phi^2 dx=N^{Q+\alpha-2}(0)\phi^2(0)=0.\] We
now obtain the desired inequality
\begin{equation}
\int _{\mathbb{G}} N^{\alpha}|\nabla_{\mathbb{G}}
N|^{\gamma}|\nabla_{\mathbb{G}} \phi|^2dx \ge
\Big(\frac{Q+\alpha-2}{2}\Big)^2  \int
_{\mathbb{G}}N^{\alpha-2}|\nabla_{\mathbb{G}} N|^{\gamma+2}\phi^2
dx.
\end{equation}

To show that the constant $(\frac{Q+\alpha-2}{2})^2$ is harp, we
use the same sequence of functions (4.6) and pass to the limit as
$\epsilon \longrightarrow 0$.
\endproof

\section{Rellich-type inequality on Carnot groups}

The classical Rellich inequality \cite{38} states that
\begin{equation} \int_{\mathbb{R}^n}|\Delta \phi(x)|^2dx\ge
\frac{n^2(n-4)^2}{16}\int_{
\mathbb{R}^n}\frac{|\phi(x)|^2}{|x|^4}dx\end{equation}
 for all
$\phi\in C_0^{\infty}(\mathbb{R}^n\setminus \{0\}$) and $n \neq
2$, where the constant $\frac{n^2(n-4)^2}{16}$   is sharp. The
Rellich inequality is the first generalization of Hardy's
inequality  to higher-order derivatives. A comprehensive study of
Rellich-type inequalities on a complete Riemannian manifold with
smooth boundary  can be found in \cite {16}. In particular, Davies
and Hinz \cite {16} obtained sharp constants $C$ for the
inequalities of the form
\[\int_{\mathbb{R}^n}\frac{|\Delta\phi(x)|^p}{|x|^{\alpha}} dx\ge
C\int_{ \mathbb{R}^n}\frac{|\phi(x)|^p}{|x|^{\beta}}dx\] for
suitable values of $\alpha, \beta, p$ and $\phi\in
C_0^{\infty}(\mathbb{R}^n\setminus \{0\}$). We should also mention
that a version of Rellich-type inequality on the Heisenberg group
has been obtained by Niu, Zhang and Wang \cite {36} and D'Ambrosio
\cite {13}. In this paper we give an analog of Rellich inequality
for general Carnot groups. The following theorem is the main
result of this section.

\begin{theorem}
Let $\mathbb{G}$ be a Carnot group with homogeneous dimension
$Q\ge 3$ and let $\phi\in C_0^{\infty}(\mathbb{G}\setminus
\{0\})$, $\alpha\in \mathbb{R}$, $Q+\alpha-4>0$. Then the
following inequality is valid
\begin{equation}
\int_{\mathbb{G}}\frac{N^{\alpha}}{|\nabla_{\mathbb{G} }N|^2}
|\Delta_{\mathbb{G}} \phi|^2dx \ge
\frac{(Q+\alpha-4)^2(Q-\alpha)^2}{16}
\int_{\mathbb{G}}N^{\alpha}\frac{|\nabla_{\mathbb{G}}
N|^2}{N^4}\phi^2dx\end{equation}
\end{theorem}
\proof A straight forward computation shows that
\begin{equation}\Delta_{\mathbb{G}}N^{\alpha-2}=(Q+\alpha-4)(\alpha-2)N^{\alpha-4}|\nabla_{\mathbb{G}}N|^2+\frac{\alpha-2}{2-Q}N^{Q+\alpha-4}\Delta
u.\end{equation} Multiplying both sides of (5.3) by $\phi^2$ and
integrating over the domain $\mathbb{G}$,  we obtain

\[\int_{\mathbb{G}}\phi^2\Delta _{\mathbb{G}}N^{\alpha-2}dx=\int_{\mathbb{G}}N^{\alpha-2}(2\phi\Delta_{\mathbb{G}}\phi+2|\nabla_{\mathbb{G}}\phi|^2)dx.\]
Since $u$ is the fundamental solution of $\Delta_{\mathbb{G}}$ and
$Q+\alpha-4>0$, we obtain
\[\int_{\mathbb{G}}\phi^2\Delta
_{\mathbb{G}}N^{\alpha-2}dx=(Q+\alpha-4)(\alpha-2)\int_{\mathbb{G}}N^{\alpha-4}|\nabla_{\mathbb{G}}N|^2\phi^2dx.\]
Therefore
\begin{equation}(Q+\alpha-4)(\alpha-2)\int_{\mathbb{G}}N^{\alpha-4}|\nabla_{\mathbb{G}}N|^2\phi^2dx-2\int_{\mathbb{G}}N^{\alpha-2}\phi\Delta_{\mathbb{G}}\phi
dx=
2\int_{\mathbb{G}}N^{\alpha-2}|\nabla_{\mathbb{G}}\phi|^2dx.\end{equation}
Applying the Hardy inequality (4.1) on the right hand side of
(5.4), we get

\[\begin{aligned}&(Q+\alpha-4)(\alpha-2)\int_{\mathbb{G}}N^{\alpha-4}|\nabla_{\mathbb{G}}N|^2\phi^2dx-2\int_{\mathbb{G}}N^{\alpha-2}\phi\Delta_{\mathbb{G}}\phi
dx\\ &\ge 2(\frac{Q+\alpha-4}{2})^2
\int_{\mathbb{G}}N^{\alpha-4}|\nabla_{\mathbb{G}}
N|^2\phi^2dx.\end{aligned}\] Now it is clear that,

\begin{equation}-\int_{\mathbb{G}}N^{\alpha-2}\phi\Delta_{\mathbb{G}}\phi dx\ge
(\frac{Q+\alpha-4}{2})(\frac{Q-\alpha}{2})\int_{\mathbb{G}}N^{\alpha-4}|\nabla_{\mathbb{G}}N|^2\phi^2dx.\end{equation}
Next, we apply the Cauchy-Schwarz inequality to the expression
$-\int_{\mathbb{G}}N^{\alpha-2}\phi\Delta\phi dx$ and we obtain
\begin{equation}-\int _{\mathbb{G}}N^{\alpha-2}\phi\Delta_{\mathbb{G}}\phi dx\le (\int
_{\mathbb{G}}N^{\alpha-4}|\nabla_{\mathbb{G}}N|^2\phi^2dx)^{1/2}(\int
_{\mathbb{G}}
\frac{|\Delta_{\mathbb{G}}\phi|^2}{|\nabla_{\mathbb{G}}N|^2}N^{\alpha}dx)^{1/2}.\end{equation}
Combining (5.5) and (5.6), we obtain the inequality (5.2).
\endproof
\noindent \textbf{Remark 4.3.} It can be shown that the constant
$C(Q,\alpha)=\frac{(Q+\alpha-4)^2(Q-\alpha)^2}{16}$  is the best
constant for the Rellich inequality (5.2), that is
\[\frac{(Q+\alpha-4)^2(Q-\alpha)^2}{16}=\inf\Big\{\frac{\int_{\mathbb{G}}N^{\alpha}\frac{|\Delta_{\mathbb{G}}f|^2}{|\nabla_{\mathbb{G}}N|^2}dx}{
\int_{\mathbb{G}}N^{\alpha}\frac{|\nabla_{\mathbb{G}}N|^2}{N^4}f^2dx},
f\in C_0^{\infty}(\mathbb{G}), f\neq 0\Big\}.\] To show this, we
modify the sequence of functions (4.6) as follows,

\[\phi_{\epsilon}(x)=\begin{cases}
 (\frac{Q+\alpha-4}{2}+\epsilon)\big(N(x)-1\big)+1 &\quad\text{if}  \quad N(x)\in [0,1],\\
N^{-(\frac{Q+\alpha-4}{2}+\epsilon)} &\quad \text{if} \quad
N(x)>1,
\end{cases}
\]
and pass to the limit as $\epsilon\longrightarrow 0$.

\noindent \textbf{Remark 4.4.} In the abelian case, when
$\mathbb{G}=\mathbb{R}^n$ with the ordinary dilations, one has
$\mathcal{G}=V_1=\mathbb{R}^n$ so that $Q=n$. Now it is clear that
the inequality (4.9)   with the homogeneous norm $N(x)=|x|$
recovers the Rellich inequality (4.8) as well as Davies-Hinz
inequality for $p=2$.
\medskip

As a consequence of Rellich-type inequality (5.2), we have the
following weighted uncertainty inequalities for sub-Laplacian
$\Delta_{\mathbb{G}}$.

\begin{corollary}Let $\mathbb{G}$ be a
Carnot group with homogeneous dimension $Q\ge 3$ and let  $\phi\in
C_0^{\infty}(\mathbb{G}\setminus \{0\})$,  $\alpha\in \mathbb{R},
Q+\alpha-4>0$. Then the following inequality valid
\begin{equation}
\Big(\int_{\mathbb{G}}
\frac{N^{4-\alpha}}{|\nabla_{\mathbb{G}}N|^2}\phi^2dx\Big)\Big(\int_{\mathbb{G}}N^{\alpha}\frac{|\Delta_{
\mathbb{G}}\phi|^2}{|\nabla_{\mathbb{G}} N|^2}dx\Big)\ge
C\Big(\int_{\mathbb{G}}\phi^2 dx\Big)^2
\end{equation}
where $C=\frac{(Q+\alpha-4)^2(Q-\alpha)^2}{16}$.
\end{corollary}

\begin{corollary}Let $\mathbb{G}$ be a
Carnot group with homogeneous dimension $Q\ge 3$ and let  $\phi\in
C_0^{\infty}(\mathbb{G}\setminus \{0\})$,  $\alpha\in \mathbb{R},
Q+\alpha-4>0$. Then the following inequality valid
\begin{equation}
\Big(\int_{\mathbb{G}}
N^{4-\alpha}|\nabla_{\mathbb{G}}N|^2\phi^2dx\Big)\Big(\int_{\mathbb{G}}N^{\alpha}\frac{|\Delta_{
\mathbb{G}}\phi|^2}{|\nabla_{\mathbb{G}} N|^2}dx\Big)\ge
C\Big(\int_{\mathbb{G}}|\nabla_{\mathbb{G}}N|^2\phi^2 dx\Big)^2
\end{equation}
where $C=\frac{(Q+\alpha-4)^2(Q-\alpha)^2}{16}$.
\end{corollary}

The following inequality and its higher order extension on the
Euclidean space  has been proved by P. P. Lions \cite {30}:
\begin{equation} \int_{\mathbb{R}^N}|\Delta \phi|^2dx \ge C\Big(
\int _{\mathbb{R}^N}\frac{\phi^q}{|x|^{2s}}dx\Big)^{2/q}, \quad
 \quad \forall\,\,\phi\in C_0^{\infty}(\mathbb{R}^N\setminus
 \{0\})\end{equation}
where, $C>0$, $N\ge 5$,  $2<q<2N/(N-4)$, and $s$ is given by
\[\frac{N-2s}{q}=\frac{N-4}{2}.\] We now obtain an analogue of the inequality (5.9)
on a metric ball in Carnot groups.
\begin{theorem} Let $\mathbb{G}$ be a
Carnot group with homogeneous dimension $Q>4$. Let $B$ be a metric
ball in $\mathbb{G}$ and $0\le s \le 2$ then for every $\phi\in
C_0^{\infty}(\mathbb{G}\setminus \{0\})$ then there exists a
constant $C>0$ such that
\[
\int_{B}\frac{|\Delta_{\mathbb{G}}\phi|^2}{|\nabla_{\mathbb{G}}N|^2}dx\ge
C\Big (\int_{B}\frac{|\nabla_{ \mathbb{G}}N|^{2s-2}}{N^{2s}}
|\phi|^{2(\frac{Q-2s}{Q-4})}dx\Big)^{\frac{Q-4}{Q-2s}}.
\]
\proof

By the H\"older inequality, we have

\begin{equation}
\begin{aligned}
\int_{B}\frac{|\nabla_{\mathbb{G}}N|^{2s-2}}{N^{2s}}\phi^{2(\frac{Q-2s}{Q-4})}dx=&
\int_{B}\frac{|\nabla_{\mathbb{G}}N|^{s}}{N^{2s}}\phi^{s}
\frac{\phi^{\frac{Q(2-s)}{Q-4}}}{|\nabla_{\mathbb{G}}
N|^{2-s}}dx\\
 \le &\Big(\int_{B}
\frac{|\nabla_{\mathbb{G}}N|^2}{N^4}\phi^2dx\Big)^{s/2}\Big(\int_{B}\frac{\phi^{\frac{2Q}{Q-4}}}{|\nabla_{\mathbb{G}}N|^2}dx\Big)^{\frac{2-s}{2}}.
\end{aligned}
\end{equation}
Using the Rellich inequality (5.2)  and weighted Sobolev
inequality \cite {34}, we get

\[\int_{B}\frac{|\nabla_{\mathbb{G}}N|^{2s-2}}{N^{2s}}\phi^{2(\frac{Q-2s}{Q-4})}dx
 \le {C_1} \Big(\int_{B}
\frac{|\Delta_{\mathbb{G}}\phi|^2}{|\nabla_{\mathbb{G}}N|^2}\Big)^{s/2}\Big(\int_{B}\frac{
|\Delta_{\mathbb{G}}\phi|^2}{|\nabla_{\mathbb{G}}N|^2}dx\Big)^{\frac{(2-s)Q}{2Q-8}}.
\]

Note that the weight function $\frac{1}{|\nabla_{\mathbb{G}}
N|^2}$ satisfies the Muckenhoupt $A_2$ condition.

Therefore

\[
\int_{B}\frac{|\Delta_{\mathbb{G}}\phi|^2}{|\nabla_{\mathbb{G}}N|^2}dx\ge
C\Big (\int_{B}\frac{|\nabla_{ \mathbb{G}}N|^{2s-2}}{N^{2s}}
|\phi|^{2(\frac{Q-2s}{Q-4})}dx\Big)^{\frac{Q-4}{Q-2s}}.
\]
\end{theorem}

\medskip

\section{Improved Hardy-type inequality} In this section we prove
Hardy-type inequalities with remainder term on Carnot groups. The
following  first theorem was inspired by the work of Brezis and
V\'azquez \cite {8} which also extends their result to Carnot
groups.

\begin{theorem} Let $\mathbb{G}$ be a Carnot group with homogeneous dimension $Q\ge 3$  and  let $B$  be a metric ball
in $ \mathbb{G}$,  $\phi\in C_0^{\infty}(\mathbb{G}\setminus
\{0\})$, $2-Q<\alpha<2$. Then the following inequality is valid
\begin{equation}\int_{\mathbb{G}} N^{\alpha}|\nabla_{\mathbb{G}} \phi|^2dx
\ge \Big(\frac{Q+\alpha-2}{2}\Big)^2 \int_{\mathbb{G}}
N^{\alpha}\frac{|\nabla_{\mathbb{G}} N|^2}{N^2}\phi
^2dx+\frac{1}{C^2 r^2(B)} \int_B\phi^2dx\end{equation}   where $C$
is a positive constant and $r(B)$ is the radius of the ball $B$.

\proof Let $\phi=N^{\frac{2-Q-\alpha}{2}}\psi$ where $\psi \in
C_0^{\infty}(B)$. We have the following result from the Theorem
4.1,

\[\int_{\mathbb{G}}N^{\alpha}|\nabla_{\mathbb{G}}\phi|^2dx
=\Big(\frac{Q+\alpha-2}{2}\Big)^2\int_{\mathbb{G}}
N^{\alpha}\frac{|\nabla_{\mathbb{G}} N|^2}{N^2}\phi ^2dx+
\int_{\mathbb{G}}N^{2-Q-\alpha}|\nabla_{\mathbb{G}}\psi|^2dx
\]
We now apply the weighted Poincar\'e inequality( see \cite{18},
\cite{32}, \cite{33})( Note that the weight function
$N^{2-Q-\alpha}$ satisfies the Muckenhoupt $A_p$ condition. See
\cite {40} for further details) and we get
\begin{equation}
\int_{\mathbb{G}}N^{\alpha}|\nabla_{\mathbb{G}}\phi|^2dx
\ge\Big(\frac{Q+\alpha-2}{2}\Big)^2\int_{\mathbb{G}}
N^{\alpha}\frac{|\nabla_{\mathbb{G}} N|^2}{N^2}\phi ^2dx+
\frac{1}{C^2 r^2(B)}\int_B N^{2-Q-\alpha}|\psi|^2dx.
\end{equation}
Therefore we have the following inequality

\begin{equation}\int_{\mathbb{G}}N^{\alpha}|\nabla_{\mathbb{G}}\phi|^2dx
\ge \Big(\frac{Q+\alpha-2}{2}\Big)^2\int_{\mathbb{G}}
N^{\alpha}\frac{|\nabla_{\mathbb{G}} N|^2}{N^2}\phi
^2dx+\frac{1}{C^2 r^2(B)}\int_B|\phi|^2 dx\end{equation} where
$1/Cr^2$ is the lower bound for the least nonzero eigenvalue of
$\Delta_{\mathbb{G}}$ on $B$.
\endproof
\end{theorem}
\medskip

The next theorem has a gradient lower order term as a remainder
term. The proof of this theorem was inspired by a recent result of
Abdellaoui, D. Colorado and I. Peral \cite {1}.

\begin{theorem}Let $\mathbb{G}$ be a Carnot group with homogeneous norm $N=u^{1/(2-Q)}$ and  let $\Omega$ be a bounded domain with smooth boundary
which contains origin, $Q\ge3$, $1<q<2$, $\phi\in
C_0^{\infty}(\Omega)$  then there exists a positive constant
$C=C(Q, q, \Omega)$ such that the following inequality is valid
\[\int_{\Omega}|\nabla_{\mathbb{G}}\phi|^2dx\ge
\Big(\frac{Q-2}{2}\Big)^2\int_{\Omega}\frac{|\nabla_{\mathbb{G}}N|^2}{N^2}\phi^2dx+C(\int_{\Omega}|\nabla_{\mathbb{G}}\phi|^qdx)^{2/q}\]
\proof Let $\psi\in C_0^{\infty}(\mathbb{G}\setminus\{0\})$ then a
straight forward computation shows that

\[|\nabla_{\mathbb{G}}\phi|^2-\nabla_{\mathbb{G}}(\frac{\phi^2}{\psi})\cdot\nabla_{\mathbb{G}}\psi=
\Big|\nabla_{\mathbb{G}}\phi-\frac{\phi}{\psi}\nabla_{\mathbb{G}}\psi\Big|^2.\]

Therefore
\[\begin{aligned}\int_{\Omega}\Big(|\nabla_{\mathbb{G}}\phi|^2-\nabla_{\mathbb{G}}(\frac{\phi^2}{\psi})\cdot\nabla_{\mathbb{G}}\psi\Big)dx &=
\int_{\Omega}\Big|\nabla_{\mathbb{G}}\phi-\frac{\phi}{\psi}\nabla_{\mathbb{G}}{\psi}\Big|^2dx\\
&\ge
c\Big(\int_{\Omega}\Big|\nabla_{\mathbb{G}}\phi-\frac{\phi}{\psi}\nabla_{\mathbb{G}}\psi\Big|^q
dx\Big)^{2/q}\end{aligned}\] where we used the Jensen's inequality
in the last step. It is clear that

\[\begin{aligned}\int_{\Omega}\Big(|\nabla_{\mathbb{G}}\phi|^2-\nabla_{\mathbb{G}}(\frac{\phi^2}{\psi})\cdot\nabla_{\mathbb{G}}\psi\Big)dx&=
\int_{\Omega}|\nabla_{\mathbb{G}}\phi|^2dx + \int_{\Omega}(\frac{\Delta_{\mathbb{G}}\psi}{\psi})\phi^2dx\\
&=\int_{\Omega}|\nabla_{\mathbb{G}}\phi|^2dx +
\beta(Q+\beta-2)\int_{\Omega}\frac{|\nabla_{\mathbb{G}}N|^2}{N^2}\phi^2dx.\end{aligned}\]
Therefore we have
\begin{equation}\int_{\Omega}|\nabla_{\mathbb{G}}\phi|^2dx\ge
-\beta(Q+\beta-2)\int_{\Omega}\frac{|\nabla_{\mathbb{G}}N|^2}{N^2}\phi^2dx+
c\Big(\int_{\Omega}\Big|\nabla_{\mathbb{G}}\phi-\frac{\phi}{\psi}\nabla_{\mathbb{G}}\psi\Big|^q
dx\Big)^{2/q}.\end{equation} Now we can use the  following
elementary inequality : Let $1<q<2$ and $w_1, w_2\in \mathbb{R}^N$
then the following inequality hold:

\begin{equation}c(q)|w_2|^q\ge |w_1+w_2|^q-|w_1|^q-q|w_1|^{q-2}\langle w_1,
w_2\rangle.\end{equation} Therefore by integration and using
successively the inequality (6.1), Young's and $L^p$-Hardy
inequalities (\cite{14}, \cite{25}) we get

\begin{equation}\int_{\Omega}\Big|\nabla_{\mathbb{G}}\phi-\frac{\phi}{\psi}\nabla_{\mathbb{G}}\psi\Big|^q
dx\ge C \int _{\Omega}|\nabla_{\mathbb{G}}\phi|^qdx.
\end{equation}
Substituting (6.6)  into (6.4)  then we get
\[\int_{\Omega}|\nabla_{\mathbb{G}}\phi|^2dx\ge -\beta(Q+\beta-2)\int_{\Omega}\frac{|\nabla_{\mathbb{G}}N|^2}{N^2}\phi^2dx+
C\Big(\int_{\Omega}|\nabla_{\mathbb{G}}\phi|^q dx\Big)^{2/q}.\]
Now choosing $\beta=\frac{2-Q}{2}$ then we obtain the desired
inequality
\[\int_{\Omega}|\nabla_{\mathbb{G}}\phi|^2dx\ge \Big(\frac{Q-2}{2}\Big)^2\int_{\Omega}\frac{|\nabla_{\mathbb{G}}N|^2}{N^2}\phi^2dx+
C\Big(\int_{\Omega}|\nabla_{\mathbb{G}}\phi|^q dx\Big)^{2/q}.\]
\end{theorem}
\endproof
\medskip

We now use the Theorem 6.2 and $L^q$  version of uncertainty
principle inequality \cite{25}, we  obtain the following
interpolation inequality.

\begin{theorem} Let $\mathbb{G}$ be a polarizable Carnot group with homogeneous norm $N=u^{1/(2-Q)}$ and  let $\Omega$ be a bounded domain with smooth boundary
which contains origin, $Q\ge3$, $1<q<2$, $1/p+1/q=1$, $C>0$. Then
for every $\phi\in C_0^{\infty}(\Omega)$  the following inequality
is valid
\[\begin{aligned}&\Big[\int_{\Omega}|\nabla_{\mathbb{G}}\phi|^2dx-
(\frac{Q-2}{2})^2\int_{\Omega}\frac{|\nabla_{\mathbb{G}}N|^2}{N^2}\phi^2dx\Big]^{1/2}\Big(\int_{\Omega}N^p|\nabla_{\mathbb{G}}N|^p\phi^pdx\Big)^{1/p}\\
&\ge\sqrt{C}\Big(\frac{Q-q}{q}\Big)\Big(\int_{\Omega}|\nabla_{\mathbb{G}}N|^2\phi^2dx\Big).\end{aligned}\]

\end{theorem}

 \noindent \textbf{Acknowledgement.} I would like to thank Jeremy Tyson for bringing to my attention the polarizable
Carnot groups and valuable discussion on these topics. I would
like to thank also Gerald Folland and Nicola Garofalo  for their
valuable comments.

\bibliographystyle{amsalpha}

\end{document}